\documentclass[11pt]{amsart}
\usepackage{amssymb, latexsym}
\theoremstyle{plain}
\newtheorem{theorem}{Theorem}

\newtheorem*{thm-cheb}{Theorem (Chebyshev)}

\newtheorem{proposition}{Proposition}

\newtheorem*{2'}{Theorem 2'}
\newtheorem*{3'}{Theorem 3'}

\theoremstyle{remark}

\newtheorem*{Remark 1}{Remark 1}
\newtheorem*{Remark 2}{Remark 2}
\newtheorem*{Remark 3}{Remark 3}
\newtheorem*{Remark 4}{Remark 4}

\numberwithin{equation}{section}

\begin{document}

\title[Secretary problem with non-uniform arrivals]
 {The secretary problem with non-uniform arrivals via a left-to-right minimum exponentially tilted distribution}

\author{Ross G. Pinsky}


\address{Department of Mathematics\\
Technion---Israel Institute of Technology\\
Haifa, 32000\\ Israel}
\email{ pinsky@math.technion.ac.il}

\urladdr{https://pinsky.net.technion.ac.il/}


\date{}

\begin{abstract}

We solve the secretary problem in the case that the ranked items arrive in a statistically biased order rather than in uniformly random order.
The bias is given by the left-to-right minimum exponentially tilted distribution with parameter $q\in(0,\infty)$.
That is, for $\sigma\in S_n$, $P_n(\sigma)$ is proportional to $q^{\text{LR}^{-}_n(\sigma)}$, where
the left-to-right minimum statistic
$\text{LR}^-_n$ is defined by
$$
\text{LR}^{-}_n(\sigma)=|\{j\in[n]: \sigma_j=\min\{\sigma_i:1\le i\le j\}\}|,\ \sigma\in S_n.
$$
For $q\in(0,1)$, higher ranked items tend to arrive earlier than in the case of the uniform distribution, and for $q\in(1,\infty)$,
  they tend to arrive later, where the highest ranked item is denoted by 1 and the lowest ranked item is denoted by $n$.
In the classical problem, the asymptotically optimal strategy is to reject the first  $M_n^*$  items, where $M_n^*\sim\frac ne$, and then to select
the first item ranked higher than any of the first  $M_n^*$ items (if such an item exists). This yields $e^{-1}$ as the  limiting   probability of success.
With the above bias on arrivals, we calculate the asymptotic behavior of the optimal strategy $M_n^*$ and the corresponding   limiting probability of success, for all regimes of $\{q_n\}_{n=1}^\infty$.
In particular, if the leading order asymptotic behavior of    $\{q_n\}_{n=1}^\infty$ is at least $\frac1{\log n}$, and if also its order is  no more than $o(n)$,
then the limiting probability of success when using an asymptotically optimal strategy  is $e^{-1}$; otherwise,
 this limiting probability of success is greater than $e^{-1}$.
Also, the limiting fraction of numbers, $\lim_{n\to\infty}\frac{M^*_n}n$, that are summarily rejected by an asymptotically optimal strategy lies in $(0,1)$ if and only if
$\lim_{n\to\infty}q_n\in(0,\infty)$.
\end{abstract}

\maketitle

\section{Introduction and Statement of  Results}
In a recent paper \cite{P22} we analyzed the secretary problem in the case that the order of arrival is biased by a Mallows distribution.
The family of Mallows distributions is obtained by exponential tilting via the inversion statistic, which introduces a bias whereby smaller  numbers tend to appear earlier
and larger numbers tend to appear later (if the parameter $q\in(0,1)$) or vice versa (if the parameter $q>1$) than in the uniform case.
In this paper we study the secretary problem with a different bias, obtained by exponential tilting via the left-to-right minimum statistic.
This latter tilting also creates a bias whereby smaller numbers tend to appear earlier and larger numbers tend to appear later (if the parameter $q\in(0,1)$) or vice versa
(if the parameter $q>1$) than in the uniform case.
It turns out that the secretary problem with bias via the left-to-right minimum statistic yields  a richer array of behavior than in the case of the Mallows distribution, and the proofs of the results require a considerably more delicate analysis
  than in the case of the Mallows distribution.

 Recall the classical secretary problem: For $n\in\mathbb{N}$, a set of $n$ ranked items is revealed, one item at a time, to an observer whose objective is to select the item with
the highest rank. The order of the items is completely random; that is, each of the  $n!$ permutations of the ranks is equally likely.
At each stage, the observer only knows the relative ranks of the items that have arrived thus far, and must  either select the current item, in which case the process terminates, or reject it and continue to the next item. If the observer rejects the first $n-1$ items, then the $n$th and final item to arrive must be accepted.
Denote by $\mathcal{S}(n,M_n)$, for $M_n\in\{0,1,\cdots, n-1\}$, the strategy whereby one rejects the first $M_n$ items and then  selects the first later arriving item
that is ranked higher than any of the first $M_n$ items (if such an item exists).
As is very well known, asymptotically as $n\to\infty$, the optimal  strategies  $\mathcal{S}(n,M_n^*)$ are those for which  $M^*_n\sim \frac ne$, and
the corresponding  limiting probability of successfully selecting the item of highest rank is $e^{-1}$.

Over the years, the secretary problem has been generalized in many directions.
 For the  secretary problem  in its classical setup, but  with items arriving in a non-uniform order, see  \cite{GM, Pf, KKN} as well as \cite{P22}.
 See \cite{GD} and \cite{GK} for some variations of the classical setup with items arriving in non-uniform order.
See \cite{B00} for a different approach to the secretary problem.
 See  \cite{F83,F89} for a history of the problem and some natural variations and generalizations.

We now define the distribution obtained by exponential tilting via the left-to-right minimum statistic.
For a permutation $\sigma\in S_n$, a number $j\in[n]$ satisfying $\sigma_j=\min\{\sigma_i:1\le i\le j\}$ is called a left-to-right minimum for $\sigma$;
note that a left-to-right minimum denotes the location of a minimum and not the value of a minimum.
The left-to-right minimum statistic
$\text{LR}^-_n$ is defined by
$$
\text{LR}^{-}_n(\sigma)=|\{j\in[n]: \sigma_j=\min\{\sigma_i:1\le i\le j\}\}|,\ \sigma\in S_n.
$$
For each $q>0$, define the left-to-right minimum exponentially tilted  distribution $P_n^{\text{LR}^-;q}$ on $S_n$ by
$$
P_n^{\text{LR}^-;q}(\sigma)= \frac{q^{\text{LR}^{-}_n(\sigma)}}{q^{(n)}},\ \sigma\in S_n,
$$
where
\begin{equation}\label{raising}
q^{(n)}:=q(q+1)\cdots(q+n-1)
\end{equation}
is the raising factorial. The fact that $q^{(n)}$ is the correct normalization constant follows from the constructions in section
\ref{construction}.

Before presenting our results   on the secretary problem, we present a simple result concerning the behavior of the expectation of the  left-to-right minimum statistic
under  $P_n^{\text{LR}^-;q_n}$ for various regimes of $\{q_n\}_{n=1}^\infty$.
\medskip

\begin{proposition}\label{lrminprop}
\noindent i. Let $q_n=o(\frac1{\log n})$.
Then
$$
\lim_{n\to\infty}E_n^{\text{LR}^-;q_n}\text{LR}^{-}_n=1.
$$
\noindent ii. Let $\lim_{n\to\infty}q_n\log n=c\in(0,\infty)$. Then
$$
\lim_{n\to\infty}E_n^{\text{LR}^-;q_n}\text{LR}^{-}_n=1+c.
$$
\noindent iii. Let $\lim_{n\to\infty}q_n\log n=\infty$ and $q_n=O(1)$. Then
$$
E_n^{\text{LR}^-;q_n}\text{LR}^{-}_n\sim q_n\log n.
$$
\noindent iv. Let $q_n\to\infty$ and $q_n=o(n)$. Then
$$
E_n^{\text{LR}^-;q_n}\text{LR}^{-}_n\sim q_n\log\frac{n+q_n}{1+q_n}.
$$
In particular, if $q_n\sim cn^\alpha$, with $c>0$ and $\alpha\in(0,1)$, then
$$
E_n^{\text{LR}^-;q_n}\text{LR}^{-}_n\sim c(1-\alpha)n^\alpha\log n.
$$
\noindent v. Let $q_n\sim cn$, with $c>0$. Then
$$
E_n^{\text{LR}^-;q_n}\text{LR}^{-}_n\sim c(\log\frac{1+c}c)n.
$$
In particular,  $c(\log\frac{1+c}c)\to\begin{cases} 0,\ \text{if}\ c\to 0;\\ 1,\ \text{if}\ c\to\infty.\end{cases}$

\noindent vi. Let $\lim_{n\to\infty}\frac{q_n}n=\infty$. Then
$$
E_n^{\text{LR}^-;q_n}\text{LR}^{-}_n\sim n.
$$
\end{proposition}
\medskip

For any permutation, the right-most location of a left-to-right minimum is the location at which the number 1 appears.
In light of this, it is intuitive from the definition of the distribution and from Proposition \ref{lrminprop} that
when $q\in(0,1)$ there is a tendency for the number 1 to appear early and when $q>1$ there is a tendency for the number 1 to appear late.
In fact,
 for $i<j$,
an exponentially tilted  distribution via the left-to-right minimum statistic has a greater effect on the placement of the number $i$ than on the placement
of the number $j$, and in particular, it has the greatest effect on  the placement of the number 1.
This tendency can be understood much
more explicitly from the first of two
constructions of $P_n^{\text{LR}^-;q_n}$ given in section \ref{construction}. In that construction, a random permutation distributed
as $P_n^{\text{LR}^-;q_n}$  is built location by location, starting with the $n$th and final location, and moving backward  one location at a time. The probability that
any number $j$ is placed in the final location is the same for all $j\in[n]-\{1\}$, but is $q$ times as much for $j=1$. Using induction, let $m\in\{1,\cdots, n-2\}$, and assume now that the locations $n,n-1,\cdots, n-m+1$ have already been filled,  say by numbers $\{i_k\}_{k=n-m+1}^n$. Then every number in $[n]-\{i_k\}_{k=n-m+1}^n$, except for the smallest one of them, has the same  probability of appearing  in location $n-m$, while the smallest of them has $q$ times as much probability to appear there.
In the final step, location 1 is filled by the one remaining number.

In light of the discussion in the above paragraph, as we turn now to the secretary  problem,
\it our convention will be that the number 1 represents the highest ranking.\rm\
Thus, for $q\in(0,1)$, there is a  tendency for the highest ranked item to arrive earlier than in the case of the uniform
distribution, while for $q>1$, their is a tendency  for it to arrive later.

If the order of arrival of the items is biased via the left-to-right minimum exponentially tilted distribution
$P_n^{\text{LR}^-;q}$
with parameter $q>0$,
let $\mathcal{P}_n^q(\mathcal{S}(n,M_n))$ denote the probability of successfully selecting the item of highest rank when employing the strategy
$\mathcal{S}(n,M_n)$, which was defined in the second paragraph of the paper.
The following theorem determines the asymptotically optimal strategies $\mathcal{S}(n,M_n^*)$ and the corresponding limiting  probability of success, for all regimes of $\{q_n\}_{n=1}^\infty$.
\medskip

\begin{theorem}\label{seclrmin}
\noindent i. Let $q_n=o(\frac1{\log n})$. Then
the asymptotically optimal strategy is $\mathcal{S}(n,M_n^*)$, where
$M_n^*=0$. (That is, the optimal strategy is to choose the first item.) The corresponding limiting probability of success is
$$
\lim_{n\to\infty}\mathcal{P}_n^q(\mathcal{S}(n,M_n^*))=1.
$$
\noindent ii. Let $q_n\sim\frac c{\log n}$, with $c\in(0,1)$.
Then
the asymptotically optimal strategy is $\mathcal{S}(n,M_n^*)$, where
$M_n^*=0$.
(That is, the optimal strategy is to choose the first item.)
The corresponding limiting probability of success is
$$
\lim_{n\to\infty}\mathcal{P}_n^q(\mathcal{S}(n,M_n^*))=e^{-c}.
$$
\noindent iii. Let $q_n\sim\frac1{\log n}$.
Then
the asymptotically optimal strategies are  $\mathcal{S}(n,M_n^*)$, where
$M_n^*=k$, for all $n\in\mathbb{N}$, where $k\in\mathbb{Z}^+$ is arbitrary,
or $\lim_{n\to\infty}M_n^*=\infty$ and $\lim_{n\to\infty}\frac{\log M_n^*}{\log n}=0$.
The
  corresponding limiting probability of success is
$$
\lim_{n\to\infty}\mathcal{P}_n^q(\mathcal{S}(n,M_n^*))=e^{-1}.
$$
\noindent iv. Let $q_n$ satisfy $lim_{n\to\infty}q_n=0$ and $\lim_{n\to\infty}q_n\log n>1$.
Then
the asymptotically optimal strategies are $\mathcal{S}(n,M_n^*)$, where
$q_n\log \frac n{M^*_n}\sim1$. (If $q_n\sim\frac c{\log  n}$ with $c>1$, then
$\lim_{n\to\infty}\frac{\log M^*_n}{\log n}=\frac{c-1}c$, and in particular, one can choose $M_n^*\sim n^{1-\frac1c}$.)
The corresponding limiting probability of success is
$$
\lim_{n\to\infty}\mathcal{P}_n^q(\mathcal{S}(n,M_n^*))=e^{-1}.
$$
\noindent v. Let    $\lim_{n\to\infty}q_n=q\in(0,\infty)$. Then the asymptotically optimal strategies are $\mathcal{S}(n,M_n^*)$, where
\begin{equation*}\label{optimali}
M_n^*\sim ne^{-\frac1q}.
\end{equation*}
 The corresponding  limiting probability of success
is
$$
\lim_{n\to\infty}\mathcal{P}_n^q(\mathcal{S}(n,M_n^*))=e^{-1}.
$$
\noindent vi. Let $q_n\to\infty$ and $q_n=o(n)$. Then
the asymptotically optimal strategies is $\mathcal{S}(n,M_n^*)$, where
\begin{equation*}\label{optimali}
n-M_n^*\sim\frac n{q_n}.
\end{equation*}
(In particular, if $q_n\sim cn^\alpha$, with $\alpha\in(0,1)$, then
$n-M_n^*\sim\frac{n^{1-\alpha}}c$.)
 The corresponding  limiting probability of success
is
$$
\lim_{n\to\infty}\mathcal{P}_n^q(\mathcal{S}(n,M_n^*))=e^{-1}.
$$
\noindent vii. Let $q_n\sim cn$, with $c\in(0,1)$.
Then
the asymptotically optimal strategy is $\mathcal{S}(n,M_n^*)$, where
\begin{equation}\label{optimalvii}
M_n^*=n-L,\ \text{if}\ \frac1L\le c<\frac1{L-1}, \ \text{where}\ 2\le L\in\mathbb{N}.
\end{equation}
The corresponding limiting probability of success is
\begin{equation}\label{optimalviiprob}
\lim_{n\to\infty}\mathcal{P}_n^q(\mathcal{S}(n,M_n^*))=\frac{cL}{(1+c)^L}, \ \text{if}\ \frac1L\le c<\frac1{L-1}, \ \text{where}\ 2\le L\in\mathbb{N}.
\end{equation}
In particular,
$$
\lim_{n\to\infty}\mathcal{P}_n^q(\mathcal{S}(n,M_n^*))>e^{-1}.
$$
\noindent viii. Let $q_n\sim cn$, with $c\ge1$.
Then
the asymptotically optimal strategy is $\mathcal{S}(n,M_n^*)$, where
$M_n^*=n-1$. (That is, the optimal strategy is to choose the last item.)
The corresponding limiting probability of success is
$$
\lim_{n\to\infty}\mathcal{P}_n^q(\mathcal{S}(n,M_n^*))=\frac c{1+c}.
$$
\noindent ix. Let $\lim_{n\to\infty}\frac{q_n}n=\infty$. Then
the asymptotically optimal strategy is $\mathcal{S}(n,M_n^*)$, where
$M_n^*=n-1$. (That is, the optimal strategy is to choose the last item.)
The corresponding limiting probability of success is
$$
\lim_{n\to\infty}\mathcal{P}_n^q(\mathcal{S}(n,M_n^*))=1.
$$
\end{theorem}
\medskip

\noindent \bf Remark 1.\rm\
 The fact that  the optimal asymptotic  probability of success is always at least  $\frac1e$ can be explained by a result of Bruss \cite{B00}.
For $n\in\mathbb{N}$, let $\{I_j\}_{j=1}^n$ be a sequence of independent indicator functions, which are observed sequentially. The observer's objective is to stop at the last $k$ for which $I_k=1$. Let $p_j$ denote the probability that $I_j=1$.
One of the results of that paper is that an optimal strategy as $n\to\infty$   yields an optimal limiting probability of at least $\frac1e$, for all choices
of $\{p_j\}_{j=1}^\infty$.
This result of Bruss can be applied to the classical secretary problem. Indeed, let $I_k$ be equal to 1 or 0 according to whether or not the $k$th item is the highest ranked item among the first $k$ items. It is easy to check that
the $\{I_k\}_{k=1}^n$ are independent under the uniform distribution. It turns out that this independence also holds under the  distributions $P_n^{\text{LR-;q}}$ (as well as under the Mallows distributions mentioned above).
The proof of this independence for $P_n^{\text{LR-;q}}$ is given in section \ref{construction}.

\medskip

\noindent \bf Remark 2.\rm\ Note that if the leading order asymptotic behavior of    $\{q_n\}_{n=1}^\infty$ is at least $\frac1{\log n}$, and if also its order is  no more than $o(n)$,
then the limiting probability of success when using an asymptotically optimal strategy  is $e^{-1}$; otherwise,
 this limiting probability of success is greater than $e^{-1}$.
Note also that the limiting fraction of numbers, $\lim_{n\to\infty}\frac{M^*_n}n$, that are summarily rejected by an asymptotically optimal strategy lies in $(0,1)$ if and only if
$\lim_{n\to\infty}q_n\in(0,\infty)$.

\medskip

\noindent\bf Remark 3.\rm\ Note the following asymmetry with respect to the cases where an optimal strategy is $M^*_n=k$, for fixed $k\in\mathbb{N}$,
and the cases where the optimal strategy is $M^*_n=n-L$, for $2\le L\in\mathbb{N}$.
For $k\in\mathbb{N}$, the strategy $M_n^*=k$ is  optimal  when $q_n\sim\frac1{\log n}$, in which case the limiting probability of success
is $e^{-1}$. However, for such $q_n$, this strategy $M_n^*=k$ is not the unique optimal strategy.
On  the other hand, for $2\le L\in\mathbb{N}$, the strategy $M_n^*=n-L$ is optimal when $q_n\sim cn$, where
$\frac1L\le c<\frac1{L-1}$. This strategy is the unique optimal strategy for such $q_n$, and the limiting probability of success
is $\frac{cL}{(1+c)^L}>e^{-1}$.
\medskip

\bf\noindent Remark 4.\rm\ As noted in the introduction, the secretary problem with bias via a Mallows distribution was analyzed in \cite{P22}.
The Mallows distributions $P_n^{\text{Mall};q}$ are obtained by exponential tilting via the inversion statistic $I_n$, which is defined by
$I_n(\sigma)=\sum_{1\le i<j\le n}1_{\sigma_j<\sigma_i}$, for $\sigma\in S_n$. Thus,
$P_n^{\text{Mall};q}(\sigma)$ is proportional to $q^{I_n(\sigma)}$. There are a variety of ways to see that
tilting via the inversion statistic has a stronger effect than tilting via the left-to-right minimum statistic.
In terms of the secretary problem, this can be seen from the fact that the limiting probability of success with left-to-right minimum tilting is $e^{-1}$
as long as $\{q_n\}_{n=1}^\infty$ behaves like $o(n)$ and is at least as large as $\frac1{\log n}$. However, as seen in \cite{P22}, for constant $q_n=q\neq1$, the limiting probability
of success is larger than $e^{-1}$.

\medskip
The following theorem gives the exact formula  for $\mathcal{P}_n^q(\mathcal{S}(n,M_n))$, for any $n,q,M_n$.
\begin{theorem}\label{lrminsecexact}
For $n\in\mathbb{N}$ and $q>0$,
\begin{equation}\label{lrminsecexactform}
\mathcal{P}_n^q(\mathcal{S}(n,M_n))=\begin{cases}q\frac{M_n}n\big(\frac{n!}{M_n!}\frac1{\prod_{l=M_n}^{n-1}(l+q)}\big)\sum_{j=M_n}^{n-1}\frac1j, \ M_n\in\{1,\cdots, n-1\};\\
\frac{(n-1)!}{\prod_{l=1}^{n-1}(l+q)},\ M_n=0.\end{cases}
\end{equation}

\end{theorem}

\medskip

The number $s(n,j)$ of permutations of $S_n$ with exactly $j$ left-to-right minima coincides with the number of permutations of $S_n$ with exactly $j$ cycles. The numbers $\{s(n,j)\}$ are called the unsigned Stirling numbers
of the first kind.   A proof of this equivalence can be given by showing that the two  quantities above satisfy the same difference equation and the same boundary conditions.
An alternative  proof is via the explicit bijection provided by Foata's  Transition Lemma \cite{B}.
This bijection maps permutations with $j$ cycles to permutations with $j$ left-to-right minima.  (Actually, using  the definition of  canonical cycle notation as presented in \cite{B},
permutations with $j$ cycles are mapped to permutations with $j$ left-to-right-maxima, but one can easily adjust the definition of  canonical cycle notation so that
permutations with $j$ cycles are mapped to permutations with $j$ left-to-right minima.)

The well-known Ewings sampling distributions are the family of distributions on $S_n$ obtained by exponential tilting via the cycle statistic.
That is, the probability of any $\sigma\in S_n$ is proportional to $q^{\text{cyc}_n(\sigma)}$, where $\text{cyc}_n(\sigma)$ denotes the number of cycles in $\sigma$.
It then follows that the distribution $P_n^{\text{LR}^-;q}$  is the push-forward distribution obtained from the Ewings sampling distribution with parameter $q$
via the bijection  from the Transition Lemma.

In order to prove
 Proposition \ref{lrminprop} and Theorem \ref{seclrmin}, it will be essential to have a so-called online construction of a random permutation
 distributed  as $P_n^{\text{LR}^-;q}$.
Such an online construction for the Ewens sampling distributions can be obtained
by a minor  tweaking  of the  classical Feller construction that builds  a uniformly random permutation  cycle by cycle
\cite{ABT, P14}.
However, combining  this construction with the push forward defined above does not  yield a useful tool for proving Proposition \ref{lrminprop} and Theorem \ref{lrminsecexact}.
 In section \ref{construction} we  give two  useful online constructions
 of a random permutation distributed according to a left-to-right minimum exponentially tilted distribution.
 The first one will be used to prove  Proposition \ref{lrminprop} and Theorem \ref{lrminsecexact}, and the second one will be used to establish the independence noted in Remark 1 after Theorem \ref{seclrmin}.

We prove
Proposition \ref{lrminprop} in section \ref{proofprop}. We prove  Theorem \ref{lrminsecexact} in section \ref{proofexact}, and then use it to prove
  Theorem \ref{seclrmin} in section \ref{proofsec}.

\section{On-line constructions of  left-to-right minimum exponentially tilted distributions}\label{construction}

We  describe  two online methods for constructing
a random permutation $\Pi^{(n)}$ distributed as $P_n^{\text{LR}^-;q}$.
Fix $q>0$.
 The first construction builds the permutation   location by location, starting with the right-most location.
For each $m\in\mathbb{N}$, define the distribution $p^{(m)}$
on $[m]$ by
\begin{equation}\label{p's}
p_i^{(m)}=\begin{cases} \frac q{q+m-1},\ i=1;\\ \frac1{q+m-1},\ i=2,\cdots, m.\end{cases}
\end{equation}
Fix $n\in\mathbb{N}$. To construct the random permutation $\Pi^{(n)}=\Pi^{(n)}_1\Pi^{(n)}_2\cdots \Pi_n^{(n)}$, distributed as $P_n^{\text{LR}^-;q}$,
make $n$ independent samples, one from each of the distributions
 $\{p^{(m)}\}_{m=1}^n$. For $m\in[n]$, denote by $\kappa_m$ the number obtained in sampling from $p^{(m)}$.
Define $\Pi^{(n)}_n=\kappa_n$. Now inductively, if $\Pi^{(n)}_n,\Pi^{(n)}_{n-1},\cdots ,\Pi^{(n)}_{m+1}$ have already been defined, let
$\Pi^{(n)}_m=\Psi_m(\kappa_m)$, where $\Psi_m$ is the increasing bijection from $[m]$ to $[n]-\{\Pi^{(n)}_k\}_{k=m+1}^n$.
Thus, for example, if $n=8$ and we sample $\kappa_8=2, \kappa_7=6,\kappa_6=1, \kappa_5=4,\kappa_4=2,\kappa_3=2,\kappa_2=1, \kappa_1=1$, then
$\Pi^{(8)}=83546172$.
By construction, the random permutation $\Pi^{(n)}$ has a left-to-right minimum at location $m$ if and only if $\kappa_m=1$.
Thus, from \eqref{p's},  for any $\sigma\in S_n$, the probability that
$\Pi^{(n)}=\sigma$ is equal to
$\frac{q^{\text{LR}^-_n(\sigma)}}{q^{(n)}}$, where
$q^{(n)}$ is as in \eqref{raising}.
\medskip

The above construction of a random permutation is a minor adaptation of the so-called $p$-shifted construction of a random permutation.
See, for example, \cite{PT} and \cite{P21}.
From Proposition 1.7 and Remark 3 following it in \cite{P21}, it follows that a $p$-shifted random permutation can also be constructed in a useful alternative fashion.
This leads to the second construction of a random permutation $\Pi^{(n)}$ with a left-to-right minimum exponentially tilted distribution.
Let $\{Y_m\}_{m=2}^\infty$ be a sequence of independent random variables with
\begin{equation}\label{Y}
P(Y_m=j)=\begin{cases}\frac q{q+m-1},\ j=0;\\ \frac1{q+m-1},\ j=1,\cdots, m-1.\end{cases}
\end{equation}
Consider now a horizontal line on which to place the numbers in $[n]$.
We begin by placing down the number 1. Then inductively, if  we have already placed down the numbers $1,2,\cdots, m-1$, the number $m$ gets placed down in the position for which there are $Y_m$ numbers to its left.
For example, for $n=8$, if $Y_2=1$, $Y_3=0$, $Y_4=1$, $Y_5=1$, $Y_6=3$, $Y_7=5$, $Y_8=7$,  then we obtain the permutation $\Pi^{(8)}=83546172$.
By the construction, for $m\in [n]$, the location of $m$ in the random permutation $\Pi^{(n)}$ will
  be a left-to-right minimum for the random permutation $\Pi^{(n)}$  if and only if $Y_m=0$. Thus, from \eqref{Y}, it follows that for any $\sigma\in S_n$, the probability that
$\Pi^{(n)}=\sigma$ is equal to $\frac{q^{\text{LR}^-_n(\sigma)}}{q^{(n)}}$.

We use this second construction now to prove the independence noted in Remark 1 after Theorem \ref{seclrmin}.
We want to
prove that for any $n\in\mathbb{N}$, the events $\{\sigma\in S_n: \sigma_m=\min(\sigma_1,\cdots, \sigma_m)\}, m=1,\cdots, n$, are independent under $P_n^{\text{LR}^-;q}$. (The event
$\{\sigma\in S_n: \sigma_m=\min(\sigma_1,\cdots, \sigma_m)\}$ is the event that $m$ is a left-to-right minimum for $\sigma$.)
It is easy to show that the number of left-to-right minima in a permutation coincides with that of its inverse; that is,
$\text{LR}^{-}_n(\sigma)=\text{LR}^{-}_n(\sigma^{-1}),\ \sigma\in S_n$.
From this fact along  with the definition of the exponentially tilted measure, it follows that if  $\sigma$ is   distributed according to  $P_n^{\text{LR}^-;q}$,
then $\sigma^{-1}$ is also  distributed according to $P_n^{\text{LR}^-;q}$. Consequently, to prove the independence of the above events under $P_n^{\text{LR}^-;q}$, it suffices to prove the independence of the events
$\{\sigma\in S_n: \sigma^{-1}_m=\min(\sigma^{-1}_1,\cdots, \sigma^{-1}_m)\}, m=1,\cdots, n$, under $P_n^{\text{LR}^-;q}$. The event
$\{\sigma\in S_n: \sigma^{-1}_m=\min(\sigma^{-1}_1,\cdots, \sigma^{-1}_m)\}$ is the event
that in the permutation $\sigma$, the number $m$ appears to the left of the numbers $1,\cdots, m-1$.
Thus,
from the second construction, this event
is the event $\{Y_m=0\}$.
This completes the proof since the $\{Y_m\}_{m=1}^n$ are independent.

\section{Proof of Proposition \ref{lrminprop}}\label{proofprop}
We use  the  first online construction in section \ref{construction} and employ the notation  from there. Under the distribution
$P_n^{\text{LR}^-;q}$,
a left-to-right minimum occurs at position
$j$ if and only  if $\kappa_j=0$,
which occurs with probability $\frac  q{j-1+q}$. Therefore
\begin{equation}\label{expectationform}
E_n^{\text{LR}^-;q_n}\text{LR}^{-}_n=\sum_{j=1}^n\frac{q_n}{j-1+q_n}=1+\sum_{j=1}^{n-1}\frac{q_n}{j+q_n}.
\end{equation}
We have
$$
\sum_{j=2}^n\frac1{j+q_n}\le\int_1^{n-1}\frac1{x+q_n}dx\le\sum_{j=1}^{n-1}\frac1{j+q_n},
$$
from which it follows that
\begin{equation}\label{intlog}
q_n\log\frac{n-1+q_n}{1+q_n}\le\sum_{j=1}^{n-1}\frac{q_n}{j+q_n}\le q_n\log\frac{n-1+q_n}{1+q_n}+\frac{q_n}{1+q_n}-\frac{q_n}{n+q_n}.
\end{equation}
Parts (i)-(v) follow almost immediately from \eqref{expectationform} and \eqref{intlog}.
Part (vi) follows from \eqref{expectationform} and \eqref{intlog} and the fact that
$\log\frac{n-1+q_n}{1+q_n}=\log(1+\frac{n-2}{1+q_n})\sim\frac n{q_n}$, for $q_n$ as in part (vi).
\hfill $\square$

\section{Proof of Theorem \ref{lrminsecexact}}\label{proofexact}
Let $\sigma=\sigma_1\sigma_2\cdots\sigma_n\in S_n$ represent the rankings of the $n$ items that arrive one by one. That is, $\sigma_j$ is the ranking of the $j$th item to arrive.
Recall that our convention is that the number 1 represents   the highest ranking.
First consider the case $M_n=0$. The strategy $\mathcal{S}(n,0)$ will select the highest ranked item if and only if $\sigma_1=1$.
We use the  first online construction in section \ref{construction}, and employ the notation from there.
The event $\{\sigma_1=1\}$  occurs if and only if $\kappa_l\neq1$, for $l=2,\cdots, n$.
Thus
$$
P_n^{\text{LR}^-;q}(\sigma_1=1)=\prod_{l=2}^n\frac{l-1}{l-1+q}.
$$
 This gives \eqref{lrminsecexactform} for the case $M_n=0$.

From now on, assume that $M_n\ge1$.
Then the strategy
$\mathcal{S}(n,M_n)$ will select the highest ranking item if and only if for some $j\in\{M_n+1,\cdots, n\}$, one has
$\sigma_j=1$ and $\min(\sigma_1,\cdots,\sigma_{j-1})=\min(\sigma_1,\cdots, \sigma_{M_n})$.
So
\begin{equation}\label{basicformula}
\mathcal{P}_n^q(\mathcal{S}(n,M_n))=\sum_{j=M_n+1}^nP_n^{\text{LR}^-;q}(\sigma_j=1,\min(\sigma_1,\cdots,\sigma_{j-1})=\min(\sigma_1,\cdots,\sigma_{M_n})).
\end{equation}
We continue to use the first  online construction in section \ref{construction}, and to employ the notation from there.
 The event $\{\sigma_j=1\}$ occurs if and only if $\kappa_l\neq1$, for $l=j+1,\cdots, n$ and $\kappa_j=1$, while the event
 $\min(\sigma_1,\cdots,\sigma_{j-1})=\min(\sigma_1,\cdots,\sigma_{M_n})$ occurs if and only  if
 $\kappa_l\neq1$, for $l=M_n+1,\cdots, j-1$. Thus,
 \begin{equation}\label{calcprob}
\begin{aligned}
& P_n^{\text{LR}^-;q}(\sigma_j=1,\min(\sigma_1,\cdots,\sigma_{j-1})=\min(\sigma_1,\cdots,\sigma_{M_n}))=\\
& \big(\prod_{l=j+1}^n\frac{l-1}{l-1+q}\big)\big(\frac q{j-1+q}\big)\big(\prod_{l=M_n+1}^{j-1}\frac{l-1}{l-1+q}\big)=\\
&\frac{q(n-1)!}{(j-1)(M_n-1)!}\frac1{\prod_{l=M_n+1}^n(l-1+q)}.
\end{aligned}
 \end{equation}
Now \eqref{lrminsecexactform} follows from \eqref{basicformula} and \eqref{calcprob}. \hfill $\square$

\section{Proof of Theorem \ref{seclrmin}}\label{proofsec}
To prove the theorem, we perform an asymptotic analysis on \eqref{lrminsecexactform} with $q=q_n$.
We begin with the estimates that are needed to  analyze the cases (v)-(ix), in which  $\{q_n\}_{n=1}^\infty$ is bounded away from zero.
Then we prove cases (v)-(ix) of the theorem. After that we prove some additional estimates that are needed for the cases (i)-(iv), in which $\lim_{n\to\infty}q_n=0$.
And then we prove cases (i)-(iv) of the theorem.

Using the well-known fact that
\begin{equation*}\label{eulergamma}
\sum_{j=1}^n\frac1j=\log n+\gamma+O(\frac1n),\ \text{where}\ \gamma\ \text{is the Euler-Mascheroni constant},
\end{equation*}
we have
\begin{equation}\label{logterm}
\sum_{j=M_n}^{n-1}\frac1j=\log \frac n{M_n}+O(\frac1{M_n}).
\end{equation}
We write
\begin{equation}\label{factorialandproduct}
\frac{n!}{M_n!}\frac1{\prod_{l=M_n}^{n-1}(l+q_n)}=\prod_{l=M_n+1}^n\frac l{l-1+q_n}
\end{equation}
Using the Taylor expansion
\begin{equation}\label{Taylor}
\log(1+x)=x-\frac12c_xx^2,\ \text{for}\ x>-1,\ \text{where}\ c_x\in(0,1),
\end{equation}
and using \eqref{logterm} for the final equality, we have
\begin{equation}\label{logprod}
\begin{aligned}
&\log\prod_{l=M_n+1}^n\frac {l-1+q_n}l=\sum_{l=M_n+1}^n\log(1+\frac{q_n-1}l)=\\
&\sum_{l=M_n+1}^n\big(\frac{q_n-1}l-\frac12c_{q_n,l}\frac{(q_n-1)^2}{l^2}\big)=\\
&(q_n-1)\log \frac n{M_n}+O(\frac{q_n-1}{M_n})+O\big((q_n-1)^2(\frac1{M_n}-\frac1n)\big),
\end{aligned}
\end{equation}
where $c_{q_n,l}\in(0,1)$.
From \eqref{factorialandproduct} and \eqref{logprod}, we have
\begin{equation}\label{productterm}
\frac{n!}{M_n!}\frac1{\prod_{l=M_n}^{n-1}(l+q_n)}=(\frac{M_n}n)^{q_n-1}\exp\Big(O(\frac{q_n-1}{M_n})+O\big((q_n-1)^2(\frac1{M_n}-\frac1n)\big)\Big).
\end{equation}
From \eqref{lrminsecexactform}, \eqref{logterm} and  \eqref{productterm}, we have
\begin{equation}\label{probanalyzed}
\begin{aligned}
&\mathcal{P}_n^{q_n}(\mathcal{S}(n,M_n))=\\
&q_n(\frac{M_n}n)^{q_n}\big(\log\frac n{M_n}+O(\frac1{M_n})\big)\exp\Big(O(\frac{q_n-1}{M_n})+O\big((q_n-1)^2(\frac1{M_n}-\frac1n)\big)\Big),\\
&\text{if}\ M_n\ge1.
\end{aligned}
\end{equation}

Using the inequality $1-x\le e^{-x}$, for $x\ge0$,
we also have
\begin{equation}\label{anotherest}
\prod_{l=M_n+1}^n\frac l{l-1+q_n}=\prod_{l=M_n+1}^n(1-\frac{q_n-1}{l-1+q_n})\le \exp\big(-(q_n-1)\sum_{l+M_n+1}^n\frac1{l-1+q_n}\big).
\end{equation}
From \eqref{logterm}, it follows that
\begin{equation}\label{harmonicest}
\sum_{l+M_n+1}^n\frac1{l-1+q_n}\ge C_{x,c}>0,\ \text{if}\  \lim_{n\to\infty}\frac{q_n}n\le c<\infty\ \text{and}\
 \lim_{n\to\infty}\frac{M_n}n\le x,
 \ \text{for}\ x\in(0,1).
\end{equation}
From \eqref{lrminsecexactform}, \eqref{logterm}, \eqref{anotherest} and \eqref{harmonicest}, we have
\begin{equation}\label{probanalyzedagain}
\begin{aligned}
&\mathcal{P}_n^{q_n}(\mathcal{S}(n,M_n))\le q_n\frac{M_n}n\big(\log\frac n{M_n}+O(\frac1{M_n})\big)\exp(-C_{x,c}(q_n-1)),
\text{where}\ C_{x,c}>0,\\
&\text{if}\ \lim_{n\to\infty}\frac{q_n}n\le c<\infty
\ \text{and}\ \lim_{n\to\infty}\frac{M_n}n\le x, \ \text{for}\ x\in(0,1).
\end{aligned}
\end{equation}

We now use the above results to  prove parts (v)-(ix).
We begin with part (v). It is easy to see that without loss of generality we can assume   that $q_n=q$ is independent of $n$.
If $\lim_{n\to\infty}\frac{M_n}n=x\in[0,1]$, then from \eqref{probanalyzed},
$$
\lim_{n\to\infty}\mathcal{P}_n^{q}(\mathcal{S}(n,M_n))=\begin{cases}-qx^q\log x, \ \text{if}\ x\in(0,1];\\ 0, \text{if}\ x=0.\end{cases}
$$
The function $-qx^q\log x$, for $x\in(0,1]$, attains its maximum value $e^{-1}$ at $x=e^{-\frac1q}$. This completes the proof of part (v).

We now prove part (vi), where we assume that $q_n\to\infty$ and $q_n=o(n)$.
It follows from \eqref{probanalyzedagain} that if $\lim_{n\to\infty}\frac{M_n}n<1$, then
$\lim_{n\to\infty}\mathcal{P}_n^{q_n}(\mathcal{S}(n,M_n))=0$. Thus,
we assume that   $\lim_{n\to\infty}\frac{M_n}n=1$ and write
\begin{equation}\label{Mn}
M_n=n-y_n, \ \text{where}\ 1\le y_n=o(n).
\end{equation}
Then from \eqref{probanalyzed}, we have
\begin{equation}\label{eston}
\begin{aligned}
&\mathcal{P}_n^{q}(\mathcal{S}(n,M_n))=q_n(1-\frac{y_n}n)^{q_n}\big(\log(1+\frac{y_n}{n-y_n})+O(\frac1n)\big)e^{o(1)}=\\
&q_n(1-\frac{y_n}n)^{q_n}\big(\frac{y_n}n+o(1)\big)e^{o(1)}.
\end{aligned}
\end{equation}
From \eqref{eston}, it follows that
\begin{equation}\label{finaleston}
\lim_{n\to\infty}\mathcal{P}_n^{q}(\mathcal{S}(n,M_n))=ze^{-z},\ \text{if}\ \lim_{n\to\infty}\frac{q_ny_n}n=z\in[0,\infty).
\end{equation}
The function $ze^{-z}$ attains its maximum value of $e^{-1}$ at $z=1$. This completes the proof of part (vi).

We now turn to  parts (vii) and (viii) together, where $q_n\sim  cn$, for some $c>0$.
In this case too  it follows from \eqref{probanalyzedagain} that if $\lim_{n\to\infty}\frac{M_n}n<1$, then
$\lim_{n\to\infty}\mathcal{P}_n^{q_n}(\mathcal{S}(n,M_n))=0$. Thus, we may  assume that $M_n$ satisfies \eqref{Mn}.
Then from \eqref{anotherest}, we have
\begin{equation}\label{anotherestagain}
\prod_{l=M_n+1}^n\frac l{l-1+q_n}\le e^{-ay_n},\ \text{for some}\ a>0.
\end{equation}
And from \eqref{lrminsecexactform}, \eqref{logterm} and \eqref{anotherestagain}, we have
\begin{equation}\label{qncnest}
\mathcal{P}_n^{q}(\mathcal{S}(n,M_n))\le q_n\frac{M_n}n\big(\log \frac n{M_n}+O(\frac1{M_n})\big)e^{-ay_n}\sim cn\big(\frac {y_n}n+o(1)\big)e^{-ay_n}.
\end{equation}
From \eqref{qncnest}, it follows that $\lim_{n\to\infty}\mathcal{P}_n^{q}(\mathcal{S}(n,M_n))=0$, if $\lim_{n\to\infty}y_n=\infty$.
Thus, we may assume now that
\begin{equation}\label{MnL}
M_n=n-L,\ L\in\mathbb{N}.
\end{equation}
From \eqref{lrminsecexactform}, we then have
\begin{equation}
\mathcal{P}_n^{q}(\mathcal{S}(n,M_n))\sim cn(1+c)^{-L}(\frac Ln)=\frac{cL}{(1+c)^L}.
\end{equation}
One has $\frac{cL}{(1+c)^L}\ge\frac{c(L+1)}{(1+c)^{L+1}}$
if and only if $c\ge\frac1L$. This
shows that if $c\in(0,1)$, then the optimal strategy is with $M_n^*$ as in \eqref{optimalvii}, and the limiting probability of success is as in
\eqref{optimalviiprob}. It also shows that if $c\ge1$, then the optimal strategy is with $M^*_n=n-1$ and the
limiting probability of success is $\frac c{1+c}$. This
 completes the proof of parts (vii) and (viii), except for the claim in (vii) that
 $\lim_{n\to\infty}\mathcal{P}_n^{q}(\mathcal{S}(n,M^*_n))>e^{-1}$.

 We now prove  this last claim. One can show that for fixed $2\le L\in\mathbb{N}$,  the expression on the right hand side of \eqref{optimalviiprob}, considered as a function
 of $c\in[\frac1L,\frac1{L-1}]$  attains its maximum value at the right hand endpoint, where it is equal to $(1+\frac1{L-1})^{-(L-1)}$.
 The claim is proved by noting that  $(1+\frac1n)^n$ increases to $e$ as $n\to\infty$.

Finally, we turn to part (ix).
The proof of this part follows from part (vi) of Proposition \ref{lrminprop}.


We now turn to the additional estimates  needed to treat the cases in which $\lim_{n\to\infty}q_n=0$.
From  \eqref{lrminsecexactform}, for $M^*_n=0$,
\begin{equation}\label{M=0}
\begin{aligned}
&\log\mathcal{P}_n^{q_n}(\mathcal{S}(n,0))=\sum_{l=1}^{n-1}\log \frac l{l+q_n}=\sum_{l=1}^{n-1}\log(1-\frac{q_n}{l+q_n})=-q_n\sum_{l=1}^{n-1}\frac1{l+q_n}+O(q_n^2)=\\
&-q_n\log n+O(q_n),\  \text{if}\ \lim_{n\to\infty}q_n=0.
\end{aligned}
\end{equation}

For fixed $M\in\mathbb{N}$, we have
\begin{equation}\label{fixedM}
\prod_{l=M+1}^n\frac {l-1+q_n}l=\frac{M+q_n}n\prod_{l=M+1}^{n-1}(1+\frac{q_n}l).
\end{equation}
Also, using \eqref{Taylor} and \eqref{logterm}, we have
\begin{equation}\label{fixedMagain}
\begin{aligned}
&\log\prod_{l=M+1}^{n-1}(1+\frac{q_n}l)=\sum_{l=M+1}^{n-1}\log(1+\frac{q_n}l)=q_n\sum_{l=M+1}^{n-1}\frac1l-\frac{q_n^2}2\sum_{l=M+1}^{n-1}\frac{c_{q_n,l}}{l^2}=\\
&q_n\big(\log n+O(1)\big)+O(q_n^2)=q_n\log n+O(q_n),\ \text{if}\ \lim_{n\to\infty}q_n=0,
\end{aligned}
\end{equation}
where $c_{q_n,l}\in(0,1)$.
From \eqref{factorialandproduct}, \eqref{fixedM} and \eqref{fixedMagain},  we have
\begin{equation}\label{Mqnto0}
\frac{n!}{M!}\frac1{\prod_{l=M}^{n-1}(l+q_n)}\sim\frac{n^{1-q_n}}M,\ \text{if}\ \lim_{n\to\infty}q_n=0,\ \text{for}\ M\in\mathbb{N}.
\end{equation}
From  \eqref{lrminsecexactform}, \eqref{logterm} and \eqref{Mqnto0}, we have
\begin{equation}\label{probMqnto0}
\begin{aligned}
&\mathcal{P}_n^{q_n}(\mathcal{S}(n,M))\sim q_n\frac Mn\frac{n^{1-q_n}}M(\log \frac nM)\sim q_n(\frac1n)^{q_n}\log n=(q_n\log n)e^{-q_n\log n},\\
& \text{if}\ \lim_{n\to\infty}q_n=0, \ \text{for}\ M\in\mathbb{N}.
\end{aligned}
\end{equation}
From \eqref{probanalyzed}, we have
\begin{equation}\label{probMnqnto0}
\begin{aligned}
&\mathcal{P}_n^{q_n}(\mathcal{S}(n,M_n))\sim q_n(\frac{M_n}n)^{q_n}\log\frac n{M_n}=(q_n\log \frac n{M_n})e^{-q_n\log \frac n{M_n}},\\
&\text{if}\ q_n\ \text{is bounded and}\ \lim_{n\to\infty}M_n=\infty.
\end{aligned}
\end{equation}

We now prove
 parts (i)-(iv). We begin with part (i), where $q_n=o(\frac1{\log n})$.
From \eqref{M=0}, \eqref{probMqnto0}
and \eqref{probMnqnto0}, and the fact that the function $xe^{-x}$ attains its maximum at $x=1$,
it follows that the optimal strategy is $\mathcal{S}(n,M_n^*)$, with $M^*_n=0$, and the
limiting probability of success is 1.
(Alternatively, part (i) follows from part (i) of Proposition \ref{lrminprop}.)

We now turn to part (ii), where $q_n\sim\frac c{\log n}$, with $c\in(0,1)$.
If we choose $M_n=M$ to be fixed, then by
\eqref{probMqnto0},
\begin{equation}\label{Mfixedc<1}
\mathcal{P}_n^{q_n}(\mathcal{S}(n,M))\sim ce^{-c}.
\end{equation}
If we choose $M_n$ such that $\lim_{n\to\infty}M_n=\infty$,
then from \eqref{probMnqnto0},
\begin{equation}\label{Mnc<1}
\mathcal{P}_n^{q_n}(\mathcal{S}(n,M))\sim c(1-\frac{\log M_n}{\log n})e^{-c(1-\frac{\log M_n}{\log n})}.
\end{equation}
If $c\in(0,1)$, the function $H_c(x)=c(1-x)e^{-c(1-x)}$ attains its maximum over $x\in[0,1]$ at $x=0$, where it is equal to $ce^{-c}$.  Thus, from \eqref{Mnc<1},
\begin{equation}\label{Mnc<1again}
\limsup_{n\to\infty}\mathcal{P}_n^{q_n}(\mathcal{S}(n,M_n))\le ce^{-c}.
\end{equation}
On the other hand,
from \eqref{M=0},
\begin{equation}\label{M=0c<1}
\lim_{n\to\infty}\mathcal{P}_n^{q_n}(\mathcal{S}(n,0))=e^{-c}.
\end{equation}
From \eqref{Mfixedc<1}, \eqref{Mnc<1again} and \eqref{M=0c<1}, if follows that the
optimal strategy is $\mathcal{S}(n,M_n^*))$, with $M^*_n=0$, and the
limiting probability of success is $e^{-c}$.

We now turn to part (iii), where $q_n\sim\frac1{\log n}$.
The analysis above for part (ii) goes through just as well when $c=1$. Thus, from the previous paragraph we conclude
that the optimal strategies $\mathcal{S}(n,M_n^*)$ are those with    $M_n^*=k\in \mathbb{Z}^+$ or
$\lim_{n\to\infty}M_n^*=\infty$ with
 $\lim_{n\to\infty}\frac{\log M_n}{\log n}=0$, and the limiting probability of success is $e^{-1}$.

We now turn to part (iv), where $\lim_{n\to\infty}q_n=0$ and $\lim_{n\to\infty}q_n\log n>1$. From \eqref{M=0}, \eqref{probMqnto0}
and \eqref{probMnqnto0}, and the fact that the function $xe^{-x}$ attains its maximum at $x=1$, it follows that
that optimal strategies are $\mathcal{S}(n,M_n^*)$, where $q_n\log \frac n{M_n}\sim1$, and the limiting probability of success is $e^{-1}$.

\hfill $\square$

\end{document}